\documentclass[reqno]{amsart}
\usepackage{amscd,amsfonts,amssymb}
\usepackage{graphicx}
\usepackage{color}

\copyrightinfo{2000}{American Mathematical Society}

\newcommand{\cD}{{\mathcal D}}

\newcommand{\cK}{{\mathcal K}}

\newcommand{\bR}{{\mathbb R}}

\numberwithin{equation}{section}

\newtheorem{Theorem}{Theorem}[section]
\newtheorem{Lemma}{Lemma}[section]

\theoremstyle{definition}
\newtheorem{Definition}{Definition}[section]
\theoremstyle{remark}
\newtheorem{Remark}{Remark}[section]

\author{V.~M.~Shelkovich}
\address{Department of Mathematics, St.-Petersburg State Architecture and
Civil Engineering University,
2 Krasnoarmeiskaya 4, 190005, St. Petersburg, Russia.}
\email{shelkv@vs1567.spb.edu}

\title[Multidimensional delta-shock waves]
{Multidimensional delta-shock waves and the transportation and concentration processes}

\thanks{The author was supported in part by DFG Project 436 RUS 113/823,
DFG Project 436 RUS 113/895, and Grant 05-01-00912 of Russian
Foundation for Basic Research.}



\date{ }

\begin{document}

\begin{abstract}
{\it $\delta$-Shock wave type solutions} in the multidimensional
system of conservation laws
$$
\rho_t + \nabla\cdot(\rho F(U))=0, \qquad
(\rho U)_t + \nabla\cdot(\rho N(U))=0,
\quad x\in \bR^n,
$$
are studied, where $F=(F_j)$ is a given vector field, $N=(N_{jk})$ is
a given tensor field, $F_j,\,N_{kj}:\bR^n \to \bR$, $j,k=1,\dots,n$;
$\rho(x,t)\in \bR$, $U(x,t)\in \bR^n$.
The well-known particular cases of this system are
zero-pressure gas dynamics in a standard form
$$
\rho_t + \nabla\cdot(\rho U)=0,
\quad
(\rho U)_t + \nabla\cdot(\rho U\otimes U)=0,
$$
and in the relativistic form
$$
\rho_t + \nabla\cdot(\rho C(U))=0,
\quad
(\rho U)_t + \nabla\cdot(\rho U\otimes C(U))=0,
$$
where $C(U)=\frac{c_0U}{\sqrt{c_0^2+|U|^2}}$,
$c_0$ is the speed of light.
We introduce the integral identities which constitute definition of $\delta$-shocks
for the above systems and using this definition derive the Rankine--Hugoniot
conditions for curvilinear $\delta$-shocks.
We show that $\delta$-shocks are connected with {\em transportation processes
and concentration processes} and derive the $\delta$-shock balance laws describing
mass and momentum transportation between the volume outside the wave
front and the wave front. In the case of zero-pressure gas dynamics
the transportation process is the concentration process.
We also prove that energy of the volume outside the wave
front and total energy are {\em nonincreasing quantities}.
The possibility of the {\em effect of kinematic self-gravitation} and
the {\em effect of dimensional bifurcations of $\delta$-shock} in
zero-pressure gas dynamics are discussed.
\end{abstract}

\maketitle

\section{Introduction}
\label{s1}

\subsection{$L^\infty$-type solutions.}\label{s1.1}
As is well known, even in the case of smooth (and, certainly,
in the case of discontinuous) initial data $U^0(x)$,
we cannot in general find a smooth solution of the one dimensional
system of conservation laws:
\begin{equation}
\label{1.1}
\left\{
\begin{array}{rclrcl}
U_t+\big(F(U)\big)_x&=&0, \quad  &&\text{in} \, \bR\times (0, \ \infty),
\medskip \\
U&=&U^0,                  \quad  &&\text{in} \, \bR\times \{t=0\}, \\
\end{array}
\right.
\end{equation}
where $F:\bR^m \to \bR^m$ is called the {\em flux-function}
associated with (\ref{1.1}); $U^0:\bR \to \bR^m$ are given smooth
vector-functions; $U=U(x,t)=(u_1(x,t),\dots,u_m(x,t))$
is the unknown function, $x\in\bR$, $t\ge 0$.

Quoting from Evans' book, ``the great difficulty in this subject is
discovering a proper notion of weak solution for the initial problem
(\ref{1.1})''~\cite[11.1.1.]{Evans}. ``We must devise some way to interpret
a less regular function $U$ as somehow ``solving'' this initial-value
problem. But as it stands, the PDE does not even make sense unless $U$
is differentiable. However, observe that if we {\em temporarily}
assume $U$ is smooth, we can as follows rewrite, so that the resulting
expression does not directly involve the derivatives of $U$. The idea
is to multiply the PDE in (\ref{1.1}) by a smooth function $\varphi$
and then to integrate by parts, thereby transferring the derivatives
onto $\varphi$''~\cite[3.4.1.a.]{Evans}.
According to the above reasoning, the following definition is introduced:
it is said that $U\in L^\infty\big(\bR\times (0,\infty);\bR^m\big)$ is a
{\em generalized solution} of the Cauchy problem (\ref{1.1})
if the integral identities
\begin{equation}
\label{2.1}
\int_{0}^{\infty}\int \Big(U\cdot{\widetilde\varphi}_{t}
+F(U)\cdot{\widetilde\varphi}_{x}\Big)\,dx\,dt
+\int U^0(x)\cdot{\widetilde\varphi}(x,0)\,dx=0
\end{equation}
hold for all compactly supported test vector-functions
${\widetilde\varphi}: \bR\times [0, \infty) \to \bR^m$,
where \ $\cdot$ \ is the scalar product of vectors, and
$\int f(x)\,dx$ denotes the improper integral
$\int_{-\infty}^{\infty}f(x)\,dx$.

Using Definition (\ref{2.1}), one can derive the classical
Rankine--Hugoniot conditions for shocks (see, e.g.,~\cite[11.1.1.]{Evans}).

\subsection{$\delta^{(n)}$-Shock wave type solutions, $n=0,1,\dots$.}\label{s1.2}
It is well known that there are ``nonclassical'' situations where, in
contrast to Lax's and Glimm's classical results, the Cauchy problem for
a system of conservation laws {\em does not possess a weak $L^{\infty}$-solution
or possesses it for some particular initial data}. In order to solve the Cauchy
problem in this ''nonclassical`` situation, it is necessary to introduce
new singular solutions called {\em $\delta$-shocks}
(see~\cite{Al-S},~\cite{B}--~\cite{Chen-Liu2},~\cite{D-S3}--~\cite{E-R-S},
~\cite{H-W},~\cite{Ke-Kr}--~\cite{Li-Y},~\cite{R1}--~\cite{Se1},
~\cite{S2}--~\cite{Y-1} and the references therein), which is a solution
such that its components {\em contain Dirac measures}.

To illustrate the above remark, we consider the Riemann problem for the
one-dimensional system of zero-pressure gas dynamics
\begin{equation}
\label{g-2-0}
\rho_t + (\rho u)_x=0, \qquad
(\rho u)_t + (\rho u^2)_x=0.
\end{equation}
If we seek a {\em shock solution} of this problem
$$
u(x,t)=u_{+}+[u]H(-x+\phi(t)),
\quad
\rho(x,t)=\rho_{+}+[\rho]H(-x+\phi(t)),
$$
we can easily verify that the velocity of the shock is
$$
\dot\phi(t)=\frac{[\rho u]}{[\rho]}\bigg|_{x=\phi(t)}
=\frac{[\rho u^2]}{[\rho u]}\bigg|_{x=\phi(t)},
$$
where $[u]=u_{-}-u_{+}$, $[\rho]=\rho_{-}-\rho_{+}$.
The last relations imply that
$$
\rho_{+}\rho_{-}\big(u_{-}-u_{+}\big)^2=0.
$$
Thus, in general, if $u_{-}\ne u_{+}$ the Riemann problem has no shock solutions.
As was written in the excellent paper by A.N.~Kraiko~\cite[page 502]{Kraiko-1},
to construct a solution for this case, it is necessary to introduce
{\em nonclassical discontinuities which may carry mass, momentum, energy}.
It is in the above-mentioned above ``nonclassical'' situations
where we have to introduce $\delta$-shocks to solve this Riemann problem
for general initial data.
Indeed, according to the above quoted papers, in the class of $\delta$-shock
type solutions this Riemann problem has a solution for general initial data
(in particular, see~\cite[Sec.~4]{D-S4}).

The theory of $\delta$-shocks has been intensively developed in the
last ten years.

Recently, in~\cite{Pan-S1},~\cite{Pan-S2},~\cite{S5}, a
concept of $\delta^{(n)}$-shock wave type solutions was introduced,
$n=1,2,\dots$. It is a {\it new type of singular solution} of a system of
conservation laws such that its components contain delta functions and
their derivatives up to $n$-th order. In~\cite{Pan-S1} the theory of
{\em $\delta'$-shocks} was established. The results~\cite{Pan-S1},
~\cite{Pan-S2},~\cite{S5} show that systems of conservation laws can
develop not only Dirac measures (as in the case of $\delta$-shocks)
but their derivatives as well.

The above singular solutions are connected with {\em transportation processes
and concentration processes}~\cite{Al-S},~\cite{Chen-Liu1},~\cite{Chen-Liu2},
~\cite{Pan-S1},~\cite{S8}.

$\delta$- and $\delta^{(n)}$-shocks, $n=1,\dots$, {do not satisfy the standard
$L^\infty$-integral identities} (\ref{2.1}). Consequently, to deal with
these singular solutions, we need

$\bullet$ to discover a {\em proper notion of a singular solution}
and to define {\em in which sense} it may satisfy a nonlinear system;

$\bullet$ to devise some way to define a {\em singular superposition
{\rm(}product{\rm)}} of distributions (for example, the product of the Heaviside
function and the delta function).

Unfortunately, using the above cited instructions from the Evans'
book~\cite[3.4.1.a.]{Evans}, $\delta^{(n)}$-shock wave type
solutions {\em cannot be defined}. Indeed, if by integrating by parts
we transfer the derivatives onto a test function $\varphi$, under the
integral sign there still remain nonlinear terms {\em undefined in the
distributional sense}, since the components of a solution may contain
Dirac measures and their derivatives.

Thus we need to develop a special technique.
Several approaches to solving $\delta$-shock problems are known
(see~\cite{Al-S},~\cite{B}--~\cite{Chen-Liu2},~\cite{D-S3}--~\cite{E-R-S},
~\cite{H-W},~\cite{Ke-Kr},~\cite{LeFl1}--~\cite{Li-Y},~\cite{R1}--~\cite{Se1},
~\cite{S2}--~\cite{S-Zh},~\cite{Y},~\cite{Y-1} and the references therein).
One of them was proposed in~\cite{D-O-S}--~\cite{D-S6}. In
these papers the {\em weak asymptotics method})
for studying the {\em dynamics of propagation and interaction}
of different singularities of quasi-linear partial differential equations
and systems of conservation laws was developed.
It appears that the {\em weak asymptotics method}
is a proper technique to deal with $\delta$- and $\delta'$-shocks.
In the framework of the {\em weak asymptotics method},
in~\cite{D-S4}--~\cite{D-S6},~\cite{S2}, definitions of $\delta$-shock wave type
solutions by integral identities for systems of conservation laws and
in~\cite{Pan-S1},~\cite{Pan-S2} the corresponding definition of $\delta'$-shock
wave type solutions were introduced. These definitions give {\em natural}
generalizations of the classical definition of the weak
$L^{\infty}$-solutions (\ref{2.1}) relevant for the structure of
$\delta$- and $\delta'$-shocks. If a solution of the Cauchy problems
contains no $\delta$ and $\delta'$-terms then these definitions
coincide with the classical definition (\ref{2.1}).
In~\cite{Al-S},~\cite{D-S3}--~\cite{D-S6},~\cite{Pan-S1},~\cite{Pan-S2},
~\cite{S2}, by using this technique, some Cauchy problems
admitting $\delta$- and $\delta'$-shocks were solved. As far as we
know, some problems related to $\delta$- and $\delta'$-shocks can
be solved only by using the {\em weak asymptotics method}.

In the numerous papers cited above $\delta$-shocks in the
{\em system of zero-pressure gas dynamics} were studied.
This is related with the fact that this system has a physical context and
is used in applications (see bellow).

In~\cite{E-R-S}, for the one-dimensional case of zero-pressure gas dynamics
(\ref{g-2-0}) a global $\delta$-shock wave type solution in the sense
of Radon measures was obtained. In~\cite{H-W}, for this system the uniqueness
of the weak solution is proved for the case when the initial value is
a Radon measure.

The multidimensional zero-pressure gas dynamics has the form
\begin{equation}
\label{g-2}
\rho_t + \nabla\cdot(\rho U)=0, \qquad
(\rho U)_t + \nabla\cdot(\rho U\otimes U)=0,
\end{equation}
where $\rho=\rho(x,t)\ge 0$ is the density, $U=(u_1(x,t),\dots,u_n(x,t))\in\bR^n$
is the velocity,
$\nabla=\big(\frac{\partial}{\partial x_1},\dots,\frac{\partial}{\partial x_n}\big)$,
\ $\cdot$ \ is the scalar product of vectors, $\otimes$ is the usual tensor
product of vectors.

In~\cite{Li-Zh1}--~\cite{Li-Y},~\cite{S-Zh},~\cite{Y-1}
the {\em planar} $\delta$-shock wave type solution in (\ref{g-2}) is defined as a
{\it measure-valued solution}. The {\it measure-valued solution} is
defined in the following way.
Let $BM(\bR^n)$ be the space of bounded Borel measures on $\bR^n$.
A pair $(\rho, U)$, where $\rho(x,t)\in C\big(BM(\bR^n),[0,\infty)\big)$,
$U(x,t)\in \big(L^\infty\big(L^\infty(\bR^n), [0,\,\infty)\big)\big)^n$,
and $U$ is measurable with respect to $\rho$ at almost all $t\ge 0$, is
said to be a {\it measure-valued solution\/} of (\ref{g-2}) in the sense
of measure if
\begin{equation}
\label{3.1*}
\begin{array}{rcl}
\displaystyle
\int_{0}^{\infty}\int_{\bR^n} \Big(\varphi_{t}+U\cdot\nabla\varphi_{x}\Big)\,d\rho\,dt&=&0, \\
\displaystyle
\int_{0}^{\infty}\int_{\bR^n} U\Big(\varphi_{t}+U\cdot\nabla\varphi_{x}\Big)\,d\rho\,dt&=&0, \\
\end{array}
\end{equation}
hold for all $\varphi(x,t) \in {\cD}(\bR^n\times[0,\,\infty))$.

In this approach a smooth discontinuity surface $\Sigma$ is parametrized as
$X=X(s)$, $t=t(s)$ ($s\in \bR^n$), separating $(X,t)$-space into two
infinite parts $\Omega_1$ and $\Omega_2$, $N=(N_X,N_t)$ is the space-time
normal to the surface $\Sigma$. The delta-shock solution takes the form
\begin{equation}
\label{3.2*-10}
\big(\rho,U\big)(X,t))=\left\{
\begin{array}{rcl}
\big(\rho_1,\,U_1\big)(X,t),            \quad  (X,t)&\in&\Omega_1, \\
\big(w(s,t)\delta(X-X(s,t),t),\,U_{\delta}(s,t)\big),  \quad (X,t)&\in&\Sigma, \\
\big(\rho_2,\,U_2\big)(X,t),            \quad  (X,t)&\in&\Omega_2. \\
\end{array}
\right.
\end{equation}
Here $U_{\delta}$ is the velocity at the points of discontinuity,
$(\rho_1,\,U_1)$ and $(\rho_2,\,U_2)$ are smooth solutions of (\ref{g-2})
in the regions $\Omega_1$ and $\Omega_2$ respectively.

In~\cite{R1},~\cite{R2}, for the 2-D case of system (\ref{g-2}) a
notion of generalized solutions in terms of the Radon measures is
introduced, and the problem of the propagation of $\delta$-shock
waves is considered. The existence of a global weak solution for
the multidimensional system of ``zero-pressure gas dynamics''
is obtained in~\cite{Se1}. The approach of the latter paper is based
on introducing of Lagrangian coordinates and on the Dafermos entropy
condition.
In~\cite{Al-S}, for a multidimensional system of zero-pressure gas
dynamics in {\it non-conservative\/} form
\begin{equation}
\label{g-1-nc}
\rho_t + \nabla\cdot(\rho U)=0,
\qquad
U_t + (U\cdot\nabla)U=0,
\end{equation}
the Cauchy problem related with propagation of a $\delta$-shock wave
was solved.
In~\cite{Danilov-1}, for multidimensional continuity equation (the first
equation in system (\ref{g-2})) the possibility of existence of
$\delta$-shock was considered.

\subsection{The physical context of zero-pressure gas dynamics.}\label{s1.3}
Study of zero-pressure
gas dynamics and its generalization is important for applications.
The zero-pressure gas dynamics can be considered as a model of the ``sticky
particle dynamics''. These models are used in many different areas of physics.
Zero-pressure gas dynamics was used to describe the formation of large-scale
structures of the universe~\cite{San-Z},~\cite{Z}; in a mathematical modeling
of pressureless mediums, in models of dusty gases (see the excellent papers
by A.N.~Kraiko and collaborators~\cite{Kraiko},~\cite{Kraiko-1}),
in modeling two-phase flows with solid particles or droplets (see the well-known
papers by A.N.~Osiptsov and collaborators~\cite{Osiptsov-1}--~\cite{Osiptsov-3}).
The presence of particles or droplets may drastically modify the flow
parameters. Moreover, large number of phenomena that are absent in
pure gas flow are inherent in two-phase flows. Among them there are
local accumulation and focusing of particles,
the inter-particle and particle-wall collisions resulting in particle
mixing and dispersion, the surface erosion due to particle impacts,
and particle-turbulence interactions which govern the dispersion
and concentration heterogeneities of inertial particles.
The dispersed phase is usually treated mathematically as a {\em pressureless
continuum}.
Zero-pressure gas dynamics was also used for modeling
the formation and evolution of traffic jams~\cite{B-D-D-R}
(F.~Berthelin, P.~Degond, M.~Delitala, M.~Rascle).

\subsection{Contents of the paper.}\label{s1.4}
In this paper we study the problems related with the $\delta$-shock in
multidimensional system of conservation laws
\begin{equation}
\label{g-1}
\rho_t + \nabla\cdot(\rho F(U))=0, \qquad
(\rho U)_t + \nabla\cdot(\rho N(U))=0,
\end{equation}
where $F=(F_1,\dots,F_n)$ is a given vector field, $N=(N_{1},\dots,N_{n})$ is
a given tensor field, $N_k=(N_{k1},\dots,N_{kn})$, $k=1,\dots,n$;
$F_j,\,N_{kj}:\bR^n \to \bR$; $\rho=\rho(x,t)$, $U=(u_1(x,t),\dots,u_n(x,t))\in\bR^n$
are the unknown functions; $x=(x_1,\dots,x_n)\in\bR^n$, $t\ge 0$.
System (\ref{g-1}) can be rewritten as
$$
\rho_t+\sum_{j=1}^n\frac{\partial }{\partial x_j}(\rho F_j(U))=0,
\quad
(\rho u_k)_t+\sum_{j=1}^n\frac{\partial}{\partial x_j}(\rho N_{kj}(U))=0,
\quad k=1,\dots,n.
$$

The well-known particular cases of this system are zero-pressure gas dynamics
in the standard form (\ref{g-2}) (here $F(U)=U$, $N(U)=U\otimes U$) and in the
relativistic form
\begin{equation}
\label{g-2-r}
\rho_t + \nabla\cdot(\rho C(U))=0,
\quad
(\rho U)_t + \nabla\cdot(\rho U\otimes C(U))=0,
\end{equation}
(here $F(U)=C(U)$, $N(U)=U\otimes C(U)$),
where $C(U)=\frac{c_0U}{\sqrt{c_0^2+|U|^2}}$, \ $c_0$ is the speed of light.
The relativistic form (\ref{g-2-r}) of zero-pressure gas dynamics was presented
in~\cite{Poupaud}.

In Sec.~\ref{s2}, we introduce the integral identities (\ref{g-4.0-10})
which constitute Definition~\ref{g-de3-1} of $\delta$-shocks for system
(\ref{g-1}). Next, using this definition, the Rankine--Hugoniot
conditions (\ref{g-51-10}) for curvilinear $\delta$-shocks are derived.

In Sec.~\ref{s3}, geometric and physical aspects of $\delta$-shocks in
system (\ref{g-1}) are studied. It is well-known that if $U\in L^\infty$
is a generalized solution of the Cauchy problem compactly supported with
respect to~$x$, then the integral $\int_{\bR^n} U(x,t)\,dx$ is independent
of time.
For $\delta$-shock wave type solutions this fact {\em does not hold}.
Nevertheless, by Theorems~\ref{g-th4} ``generalized'' analogs of these
conservation laws are derived. We prove that the ``mass'' and ``momentum''
{\em transportation processes} between the volume outside the moving
$\delta$-shock front $\Gamma_t$ and the front $\Gamma_t$ {\em are going on}.
Moreover, we derive the $\delta$-shock {\em balance relations} (\ref{g-62})
which show that the total ``mass'' $M(t)+m(t)$ and ``momentum'' $P(t)+p(t)$
are independent of time, where $M(t)$, $P(t)$ are ``mass'' and ``momentum''
of the domain outside the wave front, and $m(t)$, $p(t)$ are ``mass''
and ``momentum'' of the wave front $\Gamma_{t}$.

In Sec.~\ref{s4}, we study the case of zero-pressure gas dynamics.
The Rankine--Hugoniot conditions (\ref{g-51-11*})
for zero-pressure gas dynamics (\ref{g-2}) and the Rankine--Hugoniot
conditions (\ref{g-51-12*}) for its relativistic form (\ref{g-2-r})
are particular cases of (\ref{g-51-10}).
For zero-pressure gas dynamics system (\ref{g-2}), ``mass'' and ``momentum''
have a sense of real mass and momentum.
In this case the {\em mass transportation process} described by
Theorem~\ref{g-th4} is the {\em mass concentration process} on
the moving front $\Gamma_t$ (see Theorem~\ref{g-th5}).
According to Theorem~\ref{e-th4}, for zero-pressure gas dynamics
energy of the volume outside the wave front and the total energy are
{\em nonincreasing quantities}.
We consider the possibility of the effects of {\em kinematic self-gravitation}
and {\em dimensional bifurcations of $\delta$-shock}.

In this section we also consider a {\em spherically symmetric case of
zero-pressure gas dynamics} (\ref{g-2-sph}) and present the
Rankine--Hugoniot conditions (\ref{g-51-11*-sph}) for $\delta$-shocks
in system (\ref{g-2-sph}).
Recall that a {\em spherically symmetric case of the gas dynamics}
admits a solution which describes the {\em heavy shock}. This solution
related with the investigation of atomic bomb explosion was found by
L.~I.~Sedov, J.~von~Neumann, and G.~I.~Taylor (see~\cite[6.16.]{W}).
It seems natural that a $\delta$-shock type solution of system (\ref{g-2-sph})
can model a super explosion with a growing amplitude of the wave front.

In Appendix~\ref{s5}, some auxiliary facts are given. In particular,
we give results related with moving surfaces and distributions defined
on these surfaces, and prove the surface transport theorems.

In the author's opinion, the multidimensional system of conservation laws
(\ref{g-1}) and its generalizations can be used in physical
models which can be treated mathematically as a pressureless continuum, e.g.,
dusty gas, granular media (for example, see~\cite{Kraiko},~\cite{Kraiko-1},
~\cite{Kraiko-2},~\cite{Osiptsov-1}--~\cite{Osiptsov-3},~\cite{Aranson},
~\cite{Meerson-1},~\cite{F-Meerson2} and Subsec.~\ref{s1.3}).

\section{$\delta$-Shock type solutions and the Rankine--Hugoniot conditions}
\label{s2}

\subsection{$\delta$-Shock type solutions.}\label{s2.1}
Let $\Gamma=\bigl\{(x,t): S(x,t)=0\bigr\}$ be a hypersurface of codimension~1
in the upper half-space $\{(x,t): x\in \bR^n, \ t\in [0,\infty)\}\subset\bR^{n+1}$,
$S\in C^{\infty}(\bR^{n}\times[0,\infty))$, $\nabla S(x,t)\bigr|_{S=0}\neq 0$
for any fixed $t$, where
$\nabla=\big(\frac{\partial}{\partial x_1},\dots, \frac{\partial}{\partial x_n}\big)$.
Let $\Gamma_t=\bigl\{x: S(x,t)=0\bigr\}$ be a moving surface in $\bR^n$.
Denote by $\nu$ the unit space normal to the surface $\Gamma_t$ pointing
(in the positive direction) from $\Omega^{-}_{t}=\{x\in \bR^{n}: S(x,t)<0\}$
to $\Omega^{+}_{t}=\{x\in \bR^{n}: S(x,t)>0\}$ such that
$\nu_j=\frac{S_{x_j}}{|\nabla S|}$, \ $j=1,\dots,n$.
The direction of the vector $\nu$ coincides with the direction in which the
function $S$ increases, i.e., inward the domain $\Omega^{+}_t$.
Denote by $-G=\frac{S_t}{|\nabla S|}$ the velocity (along the normal $\nu$)
of the moving wave front $\Gamma_t$ (see Appendix~\ref{s5.1}).

For system (\ref{g-1}) we consider the {\em $\delta$-shock type initial data}
\begin{equation}
\label{5.0-10}
(U^0(x), \rho^0(x); \, U^{0}_{\delta}(x), \, x\in \Gamma_0),
\quad
\text{where}\quad \rho^0(x)={\widehat \rho}^0(x)+e^0(x)\delta(\Gamma_0),
\end{equation}
$U^0\in L^\infty\big(\bR^n;\bR^n\big)$,
${\widehat \rho}^0 \in L^\infty\big(\bR^n;\bR\big)$,
$e^0\in C(\Gamma_0)$, \ $\Gamma_0=\bigl\{x: S^0(x)=0\bigr\}$ is the
initial position of the $\delta$-shock front,
$\nabla S^0(x)\bigr|_{S^0=0}\neq 0$, \ $U^{0}_{\delta}(x)$, $x\in \Gamma_0$
is the {\em initial velocity} of the $\delta$-shock, $\delta(\Gamma_0)$
($\equiv\delta(S^0)$) is the Dirac delta function concentrated on the
surface $\Gamma_0$. The facts related to distributions defined on
surfaces can be found in Appendix~\ref{s5.2}.

Let us introduce a definition of a {\em $\delta$-shock wave
type solution} for system (\ref{g-2}).

\begin{Definition}
\label{g-de3-1} \rm
Distributions $(U,\rho)$ and a hypersurface $\Gamma$,
where $\rho(x,t)$ has the form of the sum
$$
\rho(x,t)={\widehat \rho}(x,t)+e(x,t)\delta(\Gamma),
$$
and
$U\in L^\infty\big(\bR^n\times(0,\infty);\bR^n\big)$,
${\widehat \rho} \in L^\infty\big(\bR^n\times(0,\infty);\bR\big)$,
$e\in C(\Gamma)$, is called a {\em $\delta$-shock wave type solution}
of the Cauchy problem (\ref{g-1}), (\ref{5.0-10}) if the integral identities
\begin{equation}
\label{g-4.0-10}
\begin{array}{rcl}
\displaystyle
\int_{0}^{\infty}\int
{\widehat \rho}\Big(\varphi_{t} + F(U)\cdot\nabla\varphi\Big)\,dx\,dt
+\int_{\Gamma}e \frac{\delta\varphi}{\delta t}\frac{\,d\mu(x,t)}{\sqrt{1+G^2}}
\qquad\qquad && \medskip \\
\displaystyle
+\int {\widehat\rho}^0(x)\varphi(x,0)\,dx
+\int_{\Gamma_0}e^0(x)\varphi(x,0)\,d\mu(x)&=&0, \medskip \\
\displaystyle
\int_{0}^{\infty}\int
{\widehat \rho}\Big(U\varphi_{t} + N(U)\cdot\nabla\varphi\Big)\,dx\,dt
+\int_{\Gamma}e U_{\delta}\frac{\delta\varphi}{\delta t}\frac{\,d\mu(x,t)}{\sqrt{1+G^2}}
\qquad\quad && \medskip \\
\displaystyle
+\int U^0(x){\widehat\rho}^0(x)\varphi(x,0)\,dx
+\int_{\Gamma_0}e^0(x)U_{\delta}^0(x)\varphi(x,0)\,d\mu(x)&=&0, \medskip \\
\end{array}
\end{equation}
hold for all $\varphi \in {\cD}(\bR^n\times [0, \infty))$, where
$\int f(x)\,dx$ denotes the improper integral $\int_{\bR^n}f(x)\,dx$;
\begin{equation}
\label{g-4.0-10*}
U_{\delta}=\nu G=-\frac{S_t\nabla S}{|\nabla S|^2}
\end{equation}
is the $\delta$-shock velocity, $-G=\frac{S_t}{|\nabla S|}$,
$\frac{\delta\varphi}{\delta t}$ is the $\delta$-derivative with
respect to the time variable (\ref{g-74}).
\end{Definition}

It is easy to verify that for $n=1$ Definition~\ref{g-de3-1}
coincides with the definition of $\delta$-shocks for one-dimensional zero-pressure
gas dynamics (\ref{g-2-0}) introduced in~\cite[Definition~1.2.]{D-S4}.

Let $S^0$ be a given smooth function. Denote by $\Omega^{-}_0=\{x: S^0(x)<0\}$
and $\Omega^{+}_0=\{x: S^0(x)>0\}$ the domains on the one side
and on the other side of the hypersurface $\Gamma_0=\bigl\{x: S^0(x)=0\bigr\}$.
In order to study the $\delta$-shock front-problem, i.e., to describe
the propagation of a singular front $\Gamma$ starting from the initial
position $\Gamma_0$, we need to solve the Cauchy problem for system
(\ref{g-1}) with the initial data
\begin{equation}
\label{g-3}
\begin{array}{rcl}
\displaystyle
(\rho^{0},U^0,U^{0}_{\delta}),\quad
\text{where}\quad
\rho^{0}(x)&=&\rho^{+0}(x)+[\rho^{0}(x)]H(-\Gamma_0)
+e^{0}(x)\delta(\Gamma_0), \\
\displaystyle
U^0(x)&=&U^{+0}(x)+[U^0(x)]H(-\Gamma_0). \\
\end{array}
\end{equation}
Here $[U^0(x)]=U^{-0}(x)-U^{+0}(x)$ is a jump of the function $U^0$ across
the discontinuity hypersurface $\Gamma_0$; $U^0=U^{0+}$, $\rho^0=\rho^{0+}$
if $x\in \Omega^{+}_0$, and $U^0=U^{0-}=U^{0+}+[U^0]$,
$\rho^0=\rho^{0-}=\rho^{0+}+[\rho^0]$ if $x\in \Omega^{-}_0$;
$e^{0}$ and $\rho^{0\pm}$ are given functions, $U^{0\pm}$ are given vectors;
$H(-\Gamma_0)$ ($\equiv H(-S^0)$)is the Heaviside
function defined on the surface $\Gamma_0$, $H(-\Gamma_0)=1$ if $S^0(x)<0$,
$H(-\Gamma_0)=0$ if $S^0(x)>0$.
We assume that for the initial data (\ref{g-3}) the geometric entropy condition
\begin{equation}
\label{g-4}
U^{0+}(x)\cdot \nu^0\bigr|_{\Gamma_0}
< U_{\delta}^0(x)\cdot\nu^0\bigr|_{\Gamma_0}
< U^{0-}(x)\cdot\nu^0\bigr|_{\Gamma_0}
\end{equation}
holds, where $\nu^0=\frac{\nabla S^0(x)}{|\nabla S^0(x)|}$ is the unit
space normal of $\Gamma_0$ oriented from $\Omega^{-}_0=\{x\in \bR^{n}:S^0(x)<0\}$
to $\Omega^{+}_0=\{x\in \bR^{n}: S^0(x)>0\}$.

\subsection{Rankine--Hugoniot conditions.}\label{s2.2}
Using Definition~\ref{g-de3-1}, we derive the {\em $\delta$-shock
Rankine--Hugoniot conditions} for system (\ref{g-1}).

\begin{Theorem}
\label{th1-10}
Let us assume that $\Omega\subset \bR^n\times (0, \infty)$ is a
region cut by a smooth hypersurface $\Gamma=\bigl\{(x,t): S(x,t)=0\bigr\}$
into a left- and right-hand parts $\Omega_{\mp}$. Let $(U,\rho)$, $\Gamma$
be a $\delta$-shock wave type solution of system {\rm (\ref{g-1})}
{\rm(}in the sense of Definition~{\rm\ref{g-de3-1})}, and suppose that
$(U,\rho)$ is smooth in $\Omega_{\pm}$ and has one-sided limits $U^{\pm}$,
${\widehat \rho}^{\pm}$ on $\Gamma$. Then the {\em Rankine--Hugoniot
conditions for the $\delta$-shock}
\begin{equation}
\label{g-51-10}
\begin{array}{rcl}
\displaystyle
\frac{\delta e}{\delta t}+\nabla_{\Gamma_t}\cdot(eU_{\delta})&=&
\displaystyle
\bigl([\rho F(U)],\, [\rho]\bigr)\cdot{\bf n}, \\
\displaystyle
\frac{\delta (eU_{\delta})}{\delta t}+\nabla_{\Gamma_t}\cdot(eU_{\delta}\otimes U_{\delta})&=&
\displaystyle
\bigl([\rho N(U)],\, [\rho U]\bigr)\cdot{\bf n}, \\
\end{array}
\end{equation}
hold on the discontinuity hypersurface $\Gamma$, where
${\bf n}=(\nu,-G)=\frac{\nabla_{(x,t)}S}{|\nabla S|}$
is the space-time normal to the surface $\Gamma$,
$\nabla_{(x,t)}=\big(\nabla,\frac{\partial}{\partial t}\big)$, \
$\bigl[f(U,\rho)\bigl]=f(U^{-},\rho^{-})-f(U^{+},\rho^{+})$
is a jump of the function $f(U,\rho)$ across the
discontinuity hypersurface $\Gamma$, \
$\frac{\delta}{\delta t}$ is the $\delta$-derivative {\rm(\ref{g-74})}
with respect to $t$, and the tangent gradient $\nabla_{\Gamma_t}
=\Big(\frac{\delta}{\delta x_1},\dots,\frac{\delta}{\delta x_n}\Big)$
to the surface $\Gamma_t$ is defined by {\rm(\ref{g-74})}, {\rm(\ref{g-74-3})}.
The equivalent forms of {\rm(\ref{g-51-10})} are the following:
\begin{equation}
\label{g-51-10*}
\begin{array}{rcl}
\displaystyle
\frac{\delta e}{\delta t}+\nabla_{\Gamma_t}\cdot(eU_{\delta})&=&
\displaystyle
\bigl([\rho F(U)]-[\rho]U_{\delta}\bigr)\cdot\nu, \\
\displaystyle
\frac{\delta (eU_{\delta})}{\delta t}
+\nabla_{\Gamma_t}\cdot(eU_{\delta}\otimes U_{\delta})&=&
\displaystyle
\bigl([\rho N(U)]-[\rho U]U_{\delta}\bigr)\cdot\nu, \\
\end{array}
\end{equation}
or
\begin{equation}
\label{g-51-10*-1}
\begin{array}{rcl}
\displaystyle
\frac{\delta e}{\delta t}-2{\cK} Ge&=&
\displaystyle
\bigl([\rho F(U)]-[\rho]U_{\delta}\bigr)\cdot\nu, \\
\displaystyle
\frac{\delta (eU_{\delta})}{\delta t}-2{\cK} GeU_{\delta}&=&
\displaystyle
\bigl([\rho N(U)]-[\rho U]U_{\delta}\bigr)\cdot\nu, \\
\end{array}
\end{equation}
where ${\cK}$ is the mean curvature {\rm(\ref{g-84.3})} of the moving
wave front $\Gamma_t$.
\end{Theorem}

\begin{proof}
For any test function $\varphi \in {\cD}(\Omega)$ we have
$\varphi(x,t)=0$ for $(x,t)\not\in G$, ${\overline G}\subset \Omega$.
Selecting the test function $\varphi(x,t)$ with compact support in
$\Omega_{\pm}$, we deduce from (\ref{g-4.0-10}) that (\ref{g-1}) hold in
$\Omega_{\pm}$, respectively. Now, if the test function
$\varphi(x,t)$ has the support in $\Omega$, then
$$
\int_{0}^{\infty}\int
{\widehat \rho}\Big(\varphi_{t} + F(U)\cdot\nabla\varphi\Big)\,dx\,dt
\qquad\qquad\qquad\qquad\qquad\qquad\qquad\qquad
$$
$$
=\int_{\Omega_{-}\cap G}
{\widehat \rho}\Big(\varphi_{t} + F(U)\cdot\nabla\varphi\Big)\,dx\,dt
+\int_{\Omega_{+}\cap G}
{\widehat \rho}\Big(\varphi_{t} + F(U)\cdot\nabla\varphi\Big)\,dx\,dt.
$$
Using the integrating-by-parts formula, we obtain
$$
\int_{\Omega_{\pm}\cap G}
{\widehat \rho}\Big(\varphi_{t} + F(U)\cdot\nabla\varphi\Big)\,dx\,dt
=-\int_{\Omega_{\pm}\cap G}\Big(\rho_t + \nabla\cdot(\rho F)\Big)\varphi(x,t)\,dx\,dt
$$
$$
\mp\int_{\Gamma\cap G}\Big(\frac{{\widehat \rho}^{\pm}S_t}{|\nabla_{(x,t)}S|}
+\frac{{\widehat \rho}^{\pm}F(U^{\pm})\cdot\nabla S}{|\nabla_{(x,t)}S|}\Big)
\varphi(x,t)\,d\mu(x,t)
-\int_{\Omega_{\pm}\cap G\cap \bR^n}{\widehat\rho}^0(x)\varphi(x,0)\,dx,
$$
where $d\mu(x,t)$ is the surface measure on $\Gamma$.
Next, adding the latter relations and taking into account that
$\rho_t+\nabla\cdot(\rho F)=0$, $(x,t)\in \Omega_{\pm}$, we have
$$
\int_{0}^{\infty}\int
{\widehat \rho}\Big(\varphi_{t} + F(U)\cdot\nabla\varphi\Big)\,dx\,dt
+\int {\widehat\rho}^0(x)\varphi(x,0)\,dx
\qquad\qquad
$$
\begin{equation}
\label{97-10}
\quad
=\int_{\Gamma}\Big(-[\rho]G+[\rho F(U)]\cdot\nu\Big)
\varphi(x,t)\frac{\,d\mu(x,t)}{\sqrt{1+G^2}}.
\end{equation}
Now, using the second integrating-by-parts formula in (\ref{g-84.20}), one
can see that
$$
\int_{\Gamma}e \frac{\delta\varphi}{\delta t}\frac{\,d\mu(x,t)}{\sqrt{1+G^2}}
+\int_{\Gamma_0}e^0(x)\varphi(x,0)\,d\mu(x)
=-\int_{\Gamma}\frac{\delta^* e}{\delta t}\varphi\frac{\,d\mu(x,t)}{\sqrt{1+G^2}},
$$
where the adjoint operator $\frac{\delta^* e}{\delta t}$ is defined in
(\ref{g-84.20-1}). Thus
$$
\int_{\Gamma}e \frac{\delta\varphi}{\delta t}\frac{\,d\mu(x,t)}{\sqrt{1+G^2}}
+\int_{\Gamma_0}e^0(x)\varphi(x,0)\,d\mu(x)
\qquad\qquad\qquad\qquad\qquad\qquad
$$
\begin{equation}
\label{98-10}
=-\int_{\Gamma}\Big(\frac{\delta e}{\delta t}+\nabla_{\Gamma_t}\cdot(eG\nu)\Big)
\varphi\frac{\,d\mu(x,t)}{\sqrt{1+G^2}}.
\end{equation}

Adding (\ref{97-10}) and (\ref{98-10}), we derive
$$
\int_{\Gamma}\Big(-[\rho]G+[\rho F(U)]\cdot\nu
-\frac{\delta e}{\delta t}-\nabla_{\Gamma_t}\cdot(eG\nu)\Big)
\varphi(x,t)\,\frac{d\mu(x,t)}{\sqrt{1+G^2}}=0,
$$
for all $\varphi(x,t) \in {\cD}(\Omega)$. Taking into account formula
(\ref{g-4.0-10*}) for the $\delta$-shock velocity, one can see that
the last relation implies the first relation in (\ref{g-51-10}).

In the same way as above, we obtain the second relation in (\ref{g-51-10}).

In view of (\ref{g-4.0-10*}) and (\ref{g-84.20-1}), the Rankine--Hugoniot
conditions (\ref{g-51-10}) can be rewritten as (\ref{g-51-10*}).
Since, according to the proof of Lemma~\ref{g-lem4} (see formula (\ref{g-84.20-1})),
\begin{equation}
\label{g-51-12*-kur}
\nabla_{\Gamma_t}\cdot(eU_{\delta})=-2{\cK} Ge,
\qquad
\nabla_{\Gamma_t}\cdot(eU_{\delta}\otimes U_{\delta})=-2{\cK} GeU_{\delta},
\end{equation}
the Rankine--Hugoniot conditions (\ref{g-51-10*}) can also be rewritten
in the form (\ref{g-51-10*-1}), where ${\cK}$ is the mean curvature
(\ref{g-84.3}) of the surface $\Gamma_t$.
\end{proof}

The right-hand sides of the first and second equations in
(\ref{g-51-10}) or (\ref{g-51-10*}) are called the
{\em Rankine--Hugoniot deficits} in $\rho$ and $\rho U$, respectively.

\begin{Remark}
\label{rem1-10} \rm
{\bf (a)} The Rankine--Hugoniot conditions (\ref{g-51-10}) constitute
a system of {\em second-order} PDEs. According to this fact, to solve
the Cauchy problem for system (\ref{g-1}), we use the initial data (\ref{5.0-10})
which contain the {\em initial velocity} $U^{0}_{\delta}(x)$ of a $\delta$-shock.
It is similar to the fact that in the {\em measure-valued solution}
approach~\cite{Li-W},~\cite{Li-Y},~\cite{Y-1} (see Subsec.~\ref{s1.2})
the velocity $U$ is determined on the discontinuity surface.

It remains to note that according to our Definition~\ref{g-de3-1} a
$\delta$-shock wave type solution is a pair of {\em distributions}
$(U,\rho)$ unlike the Definition~(\ref{3.1*}), where $\rho(x,t)$ is a
{\it measure\/} and $U(x,t)$ is understood as a
{\it measurable vector-function which is defined $\rho(x,t)$ a.e.\/}.

{\bf (b)} For system (\ref{g-2}) the Rankine--Hugoniot conditions
(\ref{g-51-10}) have the form (\ref{g-51-11*}). The Rankine--Hugoniot
conditions (\ref{g-51-11*}) are analogous to the Rankine--Hugoniot conditions
\begin{equation}
\label{3.3*}
\begin{array}{rcl}
\displaystyle
\frac{\partial X}{\partial t}&=&U_{\delta}(s,t), \medskip \\
\displaystyle
\frac{\partial w}{\partial t}&=&
\big([\rho U], [\rho]\big)\cdot(N_X,N_t), \medskip \\
\displaystyle
\frac{\partial (wU_{\delta})}{\partial t}&=&
\big([\rho U\otimes U], [\rho U]\big)\cdot(N_X,N_t), \\
\end{array}
\end{equation}
in the {\em measure-valued solution} approach~\cite{Li-W},
~\cite{Li-Y},~\cite{Y-1} (see Subsec.~\ref{s1.2}),
where $(N_X,N_t)$ is the space-time normal to the $\delta$-shock
front.
\end{Remark}

\section{Geometrical aspects of $\delta$-shocks: volume balance relations}
\label{s3}

It is well known that if $U\in L^\infty\big(\bR\times (0,\infty);\bR^m\big)$
is a generalized solution of the Cauchy problem (\ref{1.1}) compactly
supported with respect to~$x$, then the integral of the solution on
the whole space
\begin{equation}
\label{7.0}
\int U(x,t)\,dx=\int U^0(x)\,dx, \qquad t\ge 0
\end{equation}
is independent of time. These integrals can express the conservation
laws of {\em total area, mass, momentum, energy}, etc.
For a $\delta$-shock wave type solution the classical conservation
laws (\ref{7.0}) {\em do not hold}. However, there is a ``generalized''
analog of conservation laws (\ref{7.0}). In the one-dimensional case these
``generalized'' analogs were derived in~\cite{Al-S},~\cite{Pan-S1},~\cite{S8}.
Now we derive multidimensional generalization of these laws.

Let us assume that a moving surface $\Gamma_{t}=\bigl\{x: S(x,t)=0\bigr\}$
permanently separates $\bR^{n}_x$ into two parts $\Omega^{\pm}_{t}=\{x\in \bR^{n}: \pm S(x,t)>0\}$.
Denote $\Omega^{\pm}_{0}=\{x\in \bR^{n}: \pm S^0(x)>0\}$. Let $(U, \rho)$ be
compactly supported with respect to~$x$. Denote by
$$
M(t)=\int_{\Omega^{-}_{t}\cup\Omega^{+}_{t}}\rho(x,t)\,dx
\qquad
M(0)=\int_{\Omega^{-}_{0}\cup\Omega^{+}_{0}}\rho^0(x)\,dx,
$$
$$
P(t)=\int_{\Omega^{-}_{t}\cup\Omega^{+}_{t}}\rho(x,t)U(x,t)\,dx,
\qquad
P(0)=\int_{\Omega^{-}_{0}\cup\Omega^{+}_{0}}\rho^0(x)U^0(x)\,dx,
$$
and
$$
m(t)=\int_{\Gamma_{t}}e(x,t)\,d\mu(x), \qquad
m(0)=\int_{\Gamma_{0}}e^{0}(x)\,d\mu(x),
$$
$$
p(t)=\int_{\Gamma_{t}}e(x,t)U_{\delta}(x,t)\,d\mu(x),
\qquad
p(0)=\int_{\Gamma_{0}}e^{0}(x)U_{\delta}^0(x)\,d\mu(x),
$$
``masses'' and ``momentums'' of the domains $\Omega^{-}_{t}\cup\Omega^{+}_{t}$,
$\Omega^{-}_{0}\cup\Omega^{+}_{0}$ and the ``masses'' and ``momentums'' of the
moving wave front $\Gamma_{t}$, $\Gamma_{0}$, respectively, where $d\mu(x)$
is the surface measure on $\Gamma_{t}$. The quantities $M(t)$ and $P(t)$
can be interpreted as the volumes under the graphs $y={\widehat \rho}(x,t)$
and $Y={\widehat \rho}(x,t)U(x,t)$, \
$x\in \Omega^{-}_{t}\cup\Omega^{+}_{t}=\{x\in \bR^{n}: S(x,t)\ne 0\}$.

\begin{Theorem}
\label{g-th4}
Let $(U, \rho)$ and the discontinuity hypesurface
$\Gamma=\bigl\{(x,t): S(x,t)=0\bigr\}$ be a $\delta$-shock wave
type solution {\rm(}in the sense of Definition~{\rm\ref{g-de3-1})}
of the Cauchy problem {\rm (\ref{g-1})}, {\rm (\ref{5.0-10})}, compactly
supported with respect to~$x$,
where $\rho(x,t)={\widehat \rho}(x,t)+e(x,t)\delta(\Gamma)$.
Suppose that $(U,\rho)$ is smooth in $\Omega_{\pm}$ and has one-sided
limits $U^{\pm}$, ${\widehat \rho}^{\pm}$ on $\Gamma$.
Then the following {\em ``mass'' and ``momentum'' balance relations} hold:
\begin{equation}
\label{g-61}
\dot M(t)=-\dot m(t), \qquad
\dot P(t)=-\dot p(t);
\end{equation}
\begin{equation}
\label{g-62}
M(t)+m(t)=M(0)+m(0), \qquad P(t)+p(t)=P(0)+p(0).
\end{equation}
Thus the ``mass'' and ``momentum'' {\em transportation processes} between
the volume $\Omega^{-}_{t}\cup \Omega^{+}_{t}=\{x\in \bR^{n}: S(x,t)\ne 0\}$
and the moving front $\Gamma_t$ {\em are going on}. Moreover, the total
``mass'' $M(t)+m(t)$ and ``momentum'' $P(t)+p(t)$ are independent of time.
\end{Theorem}

\begin{proof}
Let us assume that the supports of $U(x,t)$ and $\rho(x,t)$
with respect to~$x$ belong to a compact $K\in\bR^n_x$ bounded by $\partial K$.
Let $K^{\pm}_t=\Omega^{\pm}_t\cap K$.
By $\nu$ we denote the space normal to $\Gamma_{t}$ pointing from
$\Omega^{-}_t$ to $\Omega^{+}_t$.
Differentiating $M(t)$ and using the volume transport Theorem~\ref{g-th8}, we obtain
$$
\dot M(t)=\int_{K^{-}_{t}\cup K^{+}_{t}}\frac{\partial \rho}{\partial t}\,dx
+\int_{\partial K^{-}_{t}\cup\partial K^{+}_{t}}G\rho\,d\mu(x),
$$
where $G=-\frac{S_{t}}{|\nabla S|}$.
Since $\rho^{\pm}_t + \nabla\cdot(\rho^{\pm} F(U^{\pm}))=0$,
$x\in K^{\pm}$ and the vectors $U^{\pm}$ and functions $\rho^{\pm}$
are equal to zero on the surface $\partial K^{\pm}_{t}$ except $\Gamma_{t}$,
applying Gauss's divergence theorem, we transform the last relation to the
form
$$
\dot M(t)=-\int_{K^{-}_{t}}\nabla\cdot(\rho^{-}F(U^{-}))\,dx
-\int_{K^{+}_{t}}\nabla\cdot(\rho^{+}F(U^{+}))\,dx
+\int_{\Gamma_{t}}G[\rho]\,d\mu(x)
$$
$$
\qquad\quad
=-\int_{\Gamma_{t}}\rho^{-}F(U^{-})\cdot\nu \,d\mu(x)
+\int_{\Gamma_{t}}\rho^{+}F(U^{+})\cdot\nu\ \,d\mu(x)
+\int_{\Gamma_{t}}G[\rho]\,d\mu(x)
$$
\begin{equation}
\label{g-63}
\qquad\quad
=-\int_{\Gamma_{t}}\big([\rho F(U)]\cdot\nu-[\rho]G\big)\,d\mu(x).
\end{equation}
Using the first Rankine--Hugoniot condition (\ref{g-51-10*}) and taking into
account that $G=U_{\delta}\cdot\nu$, relation (\ref{g-63}) can be rewritten as
$$
\dot M(t)=-\int_{\Gamma_{t}}\big([\rho F(U)]-[\rho]U_{\delta}\big)\cdot\nu \,d\mu(x)
=-\int_{\Gamma_{t}}\Big(\frac{\delta e}{\delta t}+\nabla_{\Gamma_t}\cdot(eU_{\delta})\Big)
\,d\mu(x).
$$
According to the surface transport Theorem~\ref{g-th4-tr}, we have
$$
\dot m(t)
=\int_{\Gamma_{t}}\Big(\frac{\delta e}{\delta t}+\nabla_{\Gamma_t}\cdot(eU_{\delta})\Big)
\,d\mu(x).
$$
Thus the first balance relation in (\ref{g-61}) is proved.

Repeating the proof of the first balance relation in (\ref{g-61})
almost word for word, we derive the second balance relation
in (\ref{g-61}).

To complete the proof of the theorem, it remains to integrate
(\ref{g-61}) with respect to~$t$.
\end{proof}

\section{Zero-pressure gas dynamics}
\label{s4}

\subsection{Rankine--Hugoniot conditions.}\label{s4.1}
According to (\ref{g-51-10}) and (\ref{g-51-10*}), for zero- pressure
gas dynamics (\ref{g-2}) and (\ref{g-2-r}) the Rankine--Hugoniot
conditions have the form
\begin{equation}
\label{g-51-11*}
\begin{array}{rcl}
\displaystyle
\frac{\delta e}{\delta t}+\nabla_{\Gamma_t}\cdot(eU_{\delta})&=&
\displaystyle
\big([\rho U]-[\rho]U_{\delta}\big)\cdot\nu, \\
\displaystyle
\frac{\delta (eU_{\delta})}{\delta t}
+\nabla_{\Gamma_t}\cdot(eU_{\delta}\otimes U_{\delta})&=&
\displaystyle
\big([\rho U\otimes U]-[\rho U]U_{\delta}\big)\cdot\nu \\
\end{array}
\end{equation}
and
\begin{equation}
\label{g-51-12*}
\begin{array}{rcl}
\displaystyle
\frac{\delta e}{\delta t}+\nabla_{\Gamma_t}\cdot(eU_{\delta})&=&
\displaystyle
\big([\rho C(U)]-[\rho]U_{\delta}\big)\cdot\nu, \\
\displaystyle
\frac{\delta (eU_{\delta})}{\delta t}
+\nabla_{\Gamma_t}\cdot(eU_{\delta}\otimes U_{\delta})&=&
\displaystyle
\big([\rho U\otimes C(U)]-[\rho U]U_{\delta}\big)\cdot\nu. \\
\end{array}
\end{equation}
respectively. Here according to (\ref{g-51-12*-kur}),
$\nabla_{\Gamma_t}\cdot(eU_{\delta})=-2{\cK} Ge$,
$\nabla_{\Gamma_t}\cdot(eU_{\delta}\otimes U_{\delta})=-2{\cK} GeU_{\delta}$.

In this case the {\em Rankine--Hugoniot deficits} in $\rho$ and $\rho U$
are the {\em currents of mass and momentum, respectively}.

{\bf Spherically symmetric case.}
It is easy to see that the solution of (\ref{g-2}) with
spherical symmetry $\rho=\rho(r,t)$, $U=u(r,t)\frac{x}{r}$, where
$r=|x|$, $x\in \bR^n$, satisfies the following system of equations
\begin{equation}
\label{g-2-sph}
\rho_t + (\rho u)_r+\frac{n-1}{r}\rho u=0, \qquad
(\rho u)_t + (\rho u^2)_r+\frac{n-1}{r}\rho u^2=0.
\end{equation}

In this case $\Gamma=\bigl\{(x,t)\in \bR^n\times[0,\infty): S(r,t)=0\bigr\}$,
$\Gamma_t=\bigl\{x\in \bR^n: S(r,t)=0\bigr\}$;
$\nabla S(r,t)=S_r\frac{x}{r}$, $|\nabla S(r,t)|=|S_r|$;
$\nu=\frac{S_r}{|S_r|}\frac{x}{r}$; $G=-\frac{S_t}{|S_r|}$;
the $\delta$-shock velocity (\ref{g-4.0-10*}) is represented as
$U_{\delta}=\nu G=-\frac{S_t}{|S_r|}\frac{x}{r}$; $x\in \bR^n$.
It is easy to see that if $f=f(r,t)$ then formulas (\ref{g-74}) have the form
\begin{equation}
\label{g-74-sph}
\frac{\delta f}{\delta t}=\frac{\partial f}{\partial t}-\frac{S_t}{S_r}
\frac{\partial f}{\partial r},
\qquad
\frac{\delta f}{\delta x_j}=0, \quad j=1,\dots,n.
\end{equation}
Now formulas (\ref{g-51-12*-kur}) read
$$
\nabla_{\Gamma_t}\cdot(eU_{\delta})=-e\frac{S_t}{S_r}\frac{n-1}{r},
\qquad
\nabla_{\Gamma_t}\cdot(eU_{\delta}\otimes U_{\delta})
=e\frac{S_t^2}{S_r^2}\frac{(n-1)x}{r^2}, \quad (x,t)\in \Gamma.
$$

Taking into account the above formula, we observe that the
Rankine--Hugoniot conditions (\ref{g-51-11*}) take the form
\begin{equation}
\label{g-51-11*-sph}
\begin{array}{rcl}
\displaystyle
e_t-\frac{S_t}{S_r}e_r-e\frac{S_t}{S_r}\frac{n-1}{r}&=&
\displaystyle
[\rho u]\frac{S_r}{|S_r|}+[\rho]\frac{S_t}{|S_r|}, \smallskip \\
\displaystyle
\Big(-e\frac{S_t}{S_r}\Big)_t-\frac{S_t}{S_r}\Big(-e\frac{S_t}{S_r}\Big)_r
+e\Big(\frac{S_t}{S_r}\Big)^2\frac{n-1}{r}&=&
\displaystyle
[\rho u^2]\frac{S_r}{|S_r|}+[\rho u]\frac{S_t}{|S_r|}, \\
\end{array}
\end{equation}
for $(x,t)\in \Gamma$.

If $S(r,t)=-r+\phi(t)$, the Rankine--Hugoniot conditions
(\ref{g-51-11*-sph}) can be rewritten as
\begin{equation}
\label{g-51-11*-sph-1}
\begin{array}{rcl}
\displaystyle
e_t+\dot\phi(t)e_r+e\dot\phi(t)\frac{n-1}{r}&=&
\displaystyle
-[\rho u]+[\rho]\dot\phi(t), \smallskip \\
\displaystyle
\Big(e\dot\phi(t)\Big)_t +\dot\phi(t)\Big(e\dot\phi(t)\Big)_r
+e\Big(\dot\phi(t)\Big)^2\frac{n-1}{r}&=&
\displaystyle
-[\rho u^2]+[\rho u]\dot\phi(t). \\
\end{array}
\end{equation}

\subsection{Physical aspects of $\delta$-shocks: mass, and momentum balance relations.}
\label{s4.3}
In this case $\rho\ge 0$ and $U$ can be considered as the gas density and
gas velocity, respectively. Here ``masses'' $M(t)$, $m(t)$ and
``momentums'' $P(t)$, $m(t)$ which were introduced in Sec.~\ref{s3}
have the sense of {\em real masses and momentums}.

To solve the Cauchy problem, we assume that for its solution
the {\em geometric entropy condition}
\begin{equation}
\label{g-55}
U^{+}(x,t)\cdot \nu\bigr|_{\Gamma_t}
< U_{\delta}(x,t)\cdot \nu\bigr|_{\Gamma_t}
< U^{-}(x,t)\cdot \nu\bigr|_{\Gamma_t},
\end{equation}
holds, where $U_{\delta}$ is the velocity (\ref{g-4.0-10*}) of the
$\delta$-shock front $\Gamma_t$, $U^{\pm}$ is the velocity behind the
$\delta$-shock wave front and ahead of it, respectively.
Condition (\ref{g-55}) implies that all characteristics on both sides
of the initial discontinuity $\Gamma_t$ must overlap. For $t=0$ the condition
(\ref{g-55}) coincides with (\ref{g-4}) .

\begin{Theorem}
\label{g-th5}
In the case of zero-pressure gas dynamics {\rm (\ref{g-2})} the
{\em transportation process} described by Theorem~{\rm\ref{g-th4}}
is the {\em mass concentration process} on the moving front $\Gamma_t$:
\begin{equation}
\label{g-61*}
\begin{array}{rcl}
\displaystyle
\dot M(t)=-\dot m(t), \,\,\,\,\quad \dot m(t)>0, &&
\qquad
\dot P(t)=-\dot p(t), \smallskip \\
\displaystyle
M(t)+m(t)=M(0)+m(0), && \qquad P(t)+p(t)=P(0)+p(0).
\end{array}
\end{equation}
\end{Theorem}

\begin{proof}
It remains to prove the inequality $\dot m(t)>0$. Since the solution
$(U, \rho)$ of the Cauchy problem (\ref{g-2}), (\ref{5.0-10}) satisfies
the entropy condition (\ref{g-55}) and $\rho^{\pm}\ge 0$, we have for
the first relation in (\ref{g-51-11*})
$$
\frac{\delta e}{\delta t}+\nabla_{\Gamma_t}\cdot(eU_{\delta})
=\big([\rho U]-[\rho]U_{\delta}\big)\cdot\nu\bigr|_{\Gamma_t}
\qquad\qquad\qquad\qquad\qquad\qquad\qquad
$$
$$
\qquad\qquad
=\big(\rho^{-}(U^{-}-U_{\delta})
\cdot\nu+\rho^{+}(U_{\delta}-U^{+})\cdot\nu\big)\bigr|_{\Gamma_t}\ge 0.
$$
This inequality and Theorem~\ref{g-th4} imply that
$\dot m(t)
=\int_{\Gamma_{t}}\big(\frac{\delta e}{\delta t}+\nabla_{\Gamma_t}\cdot(eU_{\delta})\big)
\,d\mu(x)>0$ and $\dot M(t)<0$.
In view of these inequalities, in the case of ``zero-pressure gas dynamics''
{\em mass transportation} from the volume $\Omega^{-}_{t}\cup \Omega^{+}_{t}$
to the moving wave front $\Gamma_t$ takes place.
\end{proof}

\subsection{Energy in zero-pressure gas dynamics (\ref{g-2}).}
\label{s4.4}

Let us assume that a moving surface $\Gamma_{t}=\bigl\{x: S(x,t)=0\bigr\}$
permanently separates $\bR^{n}_x$ into two parts $\Omega^{\pm}_{t}=\{x\in \bR^{n}: \pm S(x,t)>0\}$.
Let $(U, \rho)$ be compactly supported with respect to~$x$. Denote by
$$
W(t)=\frac{1}{2}\int_{\Omega^{-}_{t}\cup\Omega^{+}_{t}}\rho(x,t)|U(x,t)|^2\,dx,
\quad
w(t)=\frac{1}{2}\int_{\Gamma_{t}}e(x,t)|U_{\delta}(x,t)|^2\,d\mu(x),
$$
energies of the domain $\Omega^{-}_{t}\cup\Omega^{+}_{t}$ and of the
moving wave front $\Gamma_{t}$, respectively (see Sec.~\ref{s3}).
The function $W(t)+w(t)$ is the total energy.

\begin{Theorem}
\label{e-th4}
Let $(U, \rho)$ and the discontinuity hypesurface
$\Gamma=\bigl\{(x,t): S(x,t)=0\bigr\}$ be a $\delta$-shock wave
type solution {\rm(}in the sense of Definition~{\rm\ref{g-de3-1})}
of the Cauchy problem {\rm (\ref{g-2})}, {\rm (\ref{5.0-10})}, compactly
supported with respect to~$x$,
where $\rho(x,t)={\widehat \rho}(x,t)+e(x,t)\delta(\Gamma)$.
Suppose that $(U,\rho)$ is smooth in $\Omega_{\pm}$ and has one-sided
limits $U^{\pm}$, ${\widehat \rho}^{\pm}$ on $\Gamma$.
Then energies $W(t)$ and $W(t)+w(t)$ are nonincreasing quantities:
\begin{equation}
\label{e-61}
\frac{d}{dt}W(t)\le 0,
\qquad
\frac{d}{dt}\Big(W(t)+w(t)\Big)\le 0.
\end{equation}
\end{Theorem}

\begin{proof}
Let us assume that the supports of $U(x,t)$ and $\rho(x,t)$
with respect to~$x$ belong to a compact $K\in\bR^n_x$ bounded by $\partial K$.
Let $K^{\pm}_t=\Omega^{\pm}_t\cap K$.
By $\nu$ we denote the space normal to $\Gamma_{t}$ pointing from
$\Omega^{-}_t$ to $\Omega^{+}_t$.
Differentiating $W(t)$ and using the volume transport Theorem~\ref{g-th8}, we obtain
$$
\dot M(t)=\frac{1}{2}\Bigg(
\int_{K^{-}_{t}\cup K^{+}_{t}}\frac{\partial }{\partial t}\rho(x,t)|U(x,t)|^2\,dx
\qquad\qquad\qquad\qquad
$$
\begin{equation}
\label{e-62}
\qquad
+\int_{\partial K^{-}_{t}\cup\partial K^{+}_{t}}G\rho(x,t)|U(x,t)|^2\,d\mu(x)\bigg),
\end{equation}
where $G=-\frac{S_{t}}{|\nabla S|}$.

Since for $x\in K^{\pm}$ system (\ref{g-2}) has a smooth solution $(\rho^{\pm},U^{\pm})$,
this solution also satisfies non-conservative form (\ref{g-1-nc}). One can easily verify that
(\ref{g-2}) and (\ref{g-1-nc}) imply that
$$
(\rho^{\pm}\,|U^{\pm}|^2)_t+\nabla\cdot(\rho^{\pm}\,|U^{\pm}|^2U^{\pm})=0,
\qquad x\in K^{\pm}.
$$
Next, using the last relation, taking into account that the vectors $U^{\pm}$
and functions $\rho^{\pm}$ are equal to zero on the surface $\partial K^{\pm}_{t}$
except $\Gamma_{t}$, and applying Gauss's divergence theorem to relation
(\ref{e-62}), we transform it to the form
$$
\dot W(t)=-\int_{K^{-}_{t}}\nabla\cdot(\rho^{-}\,|U^{-}|^2U^{-})\,dx
-\int_{K^{+}_{t}}\nabla\cdot(\rho^{+}\,|U^{+}|^2U^{+})\,dx
+\int_{\Gamma_{t}}G[\rho |U|^2]\,d\mu(x)
$$
$$
\qquad\quad
=-\int_{\Gamma_{t}}\rho^{-}\,|U^{-}|^2U^{-}\cdot\nu \,d\mu(x)
+\int_{\Gamma_{t}}\rho^{+}\,|U^{+}|^2U^{+}\cdot\nu\ \,d\mu(x)
+\int_{\Gamma_{t}}G[\rho |U|^2]\,d\mu(x)
$$
\begin{equation}
\label{e-63}
\qquad\quad
=-\int_{\Gamma_{t}}\big([\rho |U|^2U]-[\rho |U|^2]U_{\delta}\big)\cdot\nu\,d\mu(x).
\end{equation}
We also take into account that $G=U_{\delta}\cdot\nu$, \ $G=-\frac{S_{t}}{|\nabla S|}$.

Since the solution $(U, \rho)$ of the Cauchy problem (\ref{g-2}), (\ref{5.0-10})
satisfies the entropy condition (\ref{g-55}) and $\rho^{\pm}\ge 0$, we have
$$
\big([\rho |U|^2U]-[\rho |U|^2]U_{\delta}\big)\cdot\nu
\qquad\qquad\qquad\qquad\qquad\qquad\qquad\qquad\qquad\qquad\qquad
$$
\begin{equation}
\label{e-63-1}
\qquad\qquad
=\big(\rho^{-}|U^{-}|^2(U^{-}-U_{\delta})\cdot\nu
+\rho^{+}|U^{+}|^2(U_{\delta}-U^{+})\cdot\nu\big)\bigr|_{\Gamma_t}\ge 0.
\end{equation}
Formulas (\ref{e-63}), (\ref{e-63-1}) imply that $\dot W(t)\le0$, i.e., the first
inequality in (\ref{e-61}) holds.

Now we will calculate $\dot w(t)$.
Taking into account formula (\ref{g-84.20-1}), due to the surface transport
Theorem~\ref{g-th4-tr}, we obtain
$$
\dot w(t)
=\frac{1}{2}\int_{\Gamma_{t}}\Big(\frac{\delta}{\delta t}\big(e(x,t)|U_{\delta}(x,t)|^2\big)
+\nabla_{\Gamma_t}\cdot(e(x,t)|U_{\delta}(x,t)|^2U_{\delta})\Big)\,d\mu(x)
$$
$$
=\frac{1}{2}\int_{\Gamma_{t}}\Big(\frac{\delta}{\delta t}\big(e(x,t)|U_{\delta}(x,t)|^2\big)
-2{\cK} Ge(x,t)|U_{\delta}(x,t)|^2\Big)\,d\mu(x)
$$
\begin{equation}
\label{e-64}
=\frac{1}{2}\int_{\Gamma_{t}}
\Big(\sum_{k=1}^n\Big(u_{\delta k}\frac{\delta(eu_{\delta k})}{\delta t}
+u_{\delta k}e\frac{\delta u_{\delta k}}{\delta t}\Big)
-2{\cK} Ge(x,t)|U_{\delta}(x,t)|^2\Big)\,d\mu(x).
\end{equation}

According to (\ref{g-51-11*}) and (\ref{g-51-12*-kur}), we have
\begin{equation}
\label{e-64-1}
\begin{array}{rcl}
\displaystyle
\frac{\delta e}{\delta t}u_{\delta k}+e\frac{\delta u_{\delta k}}{\delta t}
-2{\cK} Geu_{\delta k}&=&[\rho u_{k}U\cdot\nu]-[\rho u_{k}]U_{\delta}\cdot\nu, \\
\displaystyle
\frac{\delta e}{\delta t}u_{\delta k}
-2{\cK} Geu_{\delta k}&=&[\rho U\cdot\nu]u_{\delta k}-[\rho]U_{\delta}\cdot\nu u_{\delta k}, \\
\end{array}
\end{equation}
where $u_{\delta k}(x,t)$ is the $k$-th component of the vector $U_{\delta}$, $k=1,\dots,n$.
Now, subtracting one equation from the other in (\ref{e-64-1}), we obtain
\begin{equation}
\label{e-64-2}
e\frac{\delta u_{\delta k}}{\delta t}
=[\rho u_{k}U\cdot\nu]-[\rho u_{k}]U_{\delta}\cdot\nu
-[\rho U\cdot\nu]u_{\delta k}+[\rho]U_{\delta}\cdot\nu u_{\delta k}.
\end{equation}
Substituting equations (\ref{e-64-2}) into (\ref{e-64}), one can easily calculate
$$
\dot w(t)
=\frac{1}{2}\int_{\Gamma_{t}}
\Big(2\sum_{k=1}^n\big([\rho u_{k}U\cdot\nu]-[\rho u_{k}]U\cdot\nu\big)u_{\delta k}
\qquad\qquad\qquad\qquad
$$
$$
\qquad\qquad
-[\rho U\cdot\nu]|U_{\delta}(x,t)|^2+[\rho]|U_{\delta}(x,t)|^2U_{\delta}\cdot\nu
\Big)\,d\mu(x).
$$
Taking into account that $U_{\delta}=G\nu$, we rewrite the above relation in the form
\begin{equation}
\label{e-65}
\dot w(t)=\frac{1}{2}\int_{\Gamma_{t}}
\big(2[\rho (U\cdot\nu)^2]G-3[\rho U\cdot\nu]G^2+[\rho]G^3\big)\,d\mu(x),
\end{equation}
where $G=-\frac{S_{t}}{|\nabla S|}$.

Adding (\ref{e-63}) and (\ref{e-65}), we obtain
$$
\dot W(t)+\dot w(t)=-\frac{1}{2}\int_{\Gamma_{t}}
\big([\rho |U|^2U\cdot\nu]-[\rho |U|^2]U_{\delta}\cdot\nu
\qquad\qquad\qquad\qquad
$$
\begin{equation}
\label{e-66}
\qquad
-2[\rho (U\cdot\nu)^2]G+3[\rho U\cdot\nu]G^2-[\rho]G^3\big)\,d\mu(x).
\end{equation}
Let us represent the velocity on the wave front $U|_{\Gamma_{t}}$ as the sum
of the normal component $U\cdot\nu$ and the component $U_{tan}$ tangential
to the surface $\Gamma_{t}$. Since $|U|^2|_{\Gamma_{t}}=(U\cdot\nu)^2+U_{tan}^2$,
and $G=U_{\delta}\cdot\nu$, one can represent the integrand in (\ref{e-66}) as
$$
\big([\rho |U|^2U\cdot\nu]-[\rho |U|^2]U_{\delta}\cdot\nu
-2[\rho (U\cdot\nu)^2]G+3[\rho U\cdot\nu]G^2-[\rho]G^3\big)
$$
$$
\qquad\qquad\qquad
=\rho^{-}(U_{tan}^{-})^2(U^{-}\cdot\nu-U_{\delta}\cdot\nu)
+\rho^{+}(U_{tan}^{+})^2(U_{\delta}\cdot\nu-U^{+}\cdot\nu)
$$
\begin{equation}
\label{e-67}
\qquad\qquad\qquad
+\rho^{-}(U^{-}\cdot\nu-U_{\delta}\cdot\nu)^3
+\rho^{+}(U_{\delta}\cdot\nu-U^{-}\cdot\nu)^3.
\end{equation}
Since the solution $(U, \rho)$ of the Cauchy problem (\ref{g-2}), (\ref{5.0-10})
satisfies the entropy condition (\ref{g-55}) and $\rho^{\pm}\ge 0$, we deduce
that the last expression is {\em non-negative}.
Formulas (\ref{e-66}), (\ref{e-67}) imply that $\dot W(t)+\dot w(t)\le 0$,
i.e., the second inequality in (\ref{e-61}) holds.
\end{proof}

\begin{Remark}
\label{rem1-11} \rm
In~\cite{Chen-Liu2}, $\delta$-shock type solutions of the Riemann problem for
one- dimensional zero-pressure gas dynamics were studied as the vanishing pressure limit
of solutions to the Euler equations for nonisentropic fluids.
As stated in~\cite[p.142]{Chen-Liu2}, the limit system formally is becomes
the system of transportation equations (\ref{g-2-0}) with the additional conservation law
\begin{equation}
\label{g-2-0*}
(\rho E)_t + (\rho u E)_x=0.
\end{equation}
Next, it is stated~\cite[p.143]{Chen-Liu2} that the additional conservation law
(\ref{g-2-0*}) actually yields the entropy inequality
\begin{equation}
\label{g-2-0**}
(\rho u^2)_t + (\rho u^3)_x\le 0.
\end{equation}
in the sense of distributions for the Riemann solutions to (\ref{g-2-0}).

According to Theorem~\ref{e-th4} (for $n=1$), the inequality (\ref{g-2-0**})
{\em is not the entropy inequality}, it reflects the fact of energy nonincreasing.
\end{Remark}

\subsection{Two possible effects.}\label{s4.5}
By analyzing the $\delta$-shock Rankine--Hugoniot conditions (\ref{g-51-11*})
for zero-pressure gas dynamics (\ref{g-2}), one can deduce the possibility of
the following interesting effects.

{\bf The effect of kinematic self-gravitation.}
According to (\ref{g-61*}), in zero-pressure gas dynamics (\ref{g-2}) the mass
concentration process on the moving discontinuity surface $\Gamma_t$ is going on.
Moreover, the second $\delta$-shock Rankine--Hugoniot condition in (\ref{g-51-11*})
is the momentum conservation law.
Taking this fact into account and using Newton's second law of motion,
one can introduce an ``effective'' gravitational potential
in a neighborhood of the discontinuity surface and describe the {\em concentration
process in terms of gravitational interaction}. Since in system (\ref{g-2})
there is no term related with gravitational interaction, this ``gravitational
effect'' is of a purely {\em kinematic nature}.

{\bf Dimensional bifurcations of $\delta$-shock.}
It follows from Theorem~\ref{g-th5} that in the $n$-dimensional zero-pressure gas
dynamics (\ref{g-2}) the mass transportation process from the volume
$\Omega^{-}_{t}\cup \Omega^{+}_{t}=\{x\in \bR^{n}: S(x,t)\ne 0\}$ onto the
$n-1$-dimensional moving $\delta$-shock front $\Gamma_t$ is going on.
Let us suppose that in a {\it finite time period} $\tilde{t}$ the whole initial mass
$M(0)$ may be concentrated on $\Gamma_{t}$.
Then, according to (\ref{g-51-11*}), for $t>\tilde{t}$, {\it instead} of the
{\it whole ``initial'' $n$-dimensional} system of zero-pressure
gas dynamics (\ref{g-2}) we obtain a {\it ``surface'' $(n-1)$-dimensional
version of this system}
\begin{equation}
\label{g-51-11**}
\frac{\delta e}{\delta t}+\nabla_{\Gamma_t}\cdot(eU_{\delta})=0,
\qquad
\frac{\delta (eU_{\delta})}{\delta t}
+\nabla_{\Gamma_t}\cdot(eU_{\delta}\otimes U_{\delta})=0,
\end{equation}
where instead of the gas velocity $U$ we have the velocity $U_{\delta}$
of the moving $\delta$-shock front $\Gamma_t$, and instead of the gas volume
density $\rho$ we have the {\it surface density of the front mass} $e$.
Moreover, the quantities $U_{\delta}$, $e$ are defined only on the moving
front $\Gamma_t$. Since system (\ref{g-51-11**}) is an $n-1$-dimensional
analog of system (\ref{g-2}) on the $(n-1)$-- dimensional surface $\Gamma_{t}$
as on a Riemannian manifold, therefore its solution can develop singularities
within a finite time period, and the whole mass is concentrated on the manifold of
dimension $n-2$, and so on. Thus, it may happen that after a finite
number of bifurcations the whole initial mass  will be concentrated at a
singular point.

A description of the above effect prompts to generalize Definition~\ref{g-de3-1}
and introduce a new concept of a multidimensional $\delta$-shock type solution to
system (\ref{g-1}) as a pair $(U,\rho)$, where $\rho(x,t)$ has the form of the sum
\begin{equation}
\label{g-4.0-10***}
\rho(x,t)={\widehat \rho}(x,t)+\sum_{j=1}^ne_j(x,t)\delta(\Gamma^{(j)}),
\end{equation}
$U\in L^\infty\big(\bR^n\times(0,\infty);\bR^n\big)$,
${\widehat \rho} \in L^\infty\big(\bR^n\times(0,\infty);\bR\big)$,
$e_j\in C(\Gamma^{(j)})$, \ $\Gamma^{(j)}$ is a hypersurface of
codimension $j$, \ $\delta(\Gamma^{(j)})$ the Dirac delta function
concentrated on the hypersurface $\Gamma^{(j)}$, \ $j=1,2,\dots,n$.
For this purpose we need to derive special integral identities
analogous to (\ref{g-4.0-10}) and develop the theory of such type of solutions.
In the framework of such type definition one can solve the above problem
of dimensional bifurcations of $\delta$-shock.

\appendix

\section{Some auxiliary facts}
\label{s5}

\subsection{Moving surfaces of discontinuity.}\label{s5.1}
Let us present some results concerning moving surfaces
from~\cite[5.2.]{Kan},~\cite{Be1},~\cite{Be2}.
Let $\Gamma_t$ be a smooth moving surface of codimension~1 in the
space $\bR^n$. Such a surface can be represented locally either
in the form $\Gamma_t=\bigl\{x\in \bR^n: S(x,t)=0\bigr\}$, or in
terms of the curvilinear Gaussian coordinates $s=(s_1,\dots,s_{n-1})$
on the surface:
$$
x_j=x_j(s_1,\dots,s_{n-1},t), \qquad s\in \bR^{n-1}.
$$
We also consider the surface $\Gamma=\bigl\{(x,t)\in \bR^{n+1}: S(x,t)=0\bigr\}$
as a submanifold of the space-time $\bR^n\times\bR$.
We shall assume that $\nabla S(x,t)\bigr|_{\Gamma_t}\neq 0$ for all fixed values
of $t$, where $\nabla
=\big(\frac{\partial}{\partial x_1},\dots,\frac{\partial}{\partial x_n}\big)$.
Let $\nu$ be the unit space normal to the surface $\Gamma_t$ pointing in
the positive direction such that
$\frac{\partial S}{\partial x_j}=|\nabla S|\nu_j$, \ $j=1,\dots,n$.

Let $f(x,t)$ be a function defined on the surface $\Gamma_t$
for some time interval, and denote by $\frac{\delta f}{\delta t}$
the derivative with respect to time as it would be computed by an
observer moving with the surface. This derivative has the following
geometrical interpretation. Let $M_0$ be a point on the surface at
the time $t=t_0$. Construct the normal line to the surface at $M_0$.
At the time $t=t_0+\Delta t$, $\Delta t$ is an infinitesimal, this normal
meets the surface $\Gamma_{t+\Delta t}$ at the point $M=M(t+\Delta t)$.
Then the $\delta$-derivative is defined as
\begin{equation}
\label{g-70}
\frac{\delta f(M_0,t_0)}{\delta t}=\lim_{\Delta t \to 0}
\frac{f(M)-f(M_0)}{\Delta t}.
\end{equation}
If $\Delta s$ is the distance between $M_0$ and $M$, then
\begin{equation}
\label{g-71}
G=\lim_{\Delta t \to 0}\frac{\Delta s}{\Delta t}
\end{equation}
is the {\em normal velocity of the moving surface} $\Gamma_t$ and
\begin{equation}
\label{g-72}
\frac{\delta x_j}{\delta t}=\lim_{\Delta t \to 0}\frac{\Delta x_j}{\Delta t}
=\lim_{\Delta t \to 0}\frac{\Delta s}{\Delta t}\frac{\Delta x_j}{\Delta s}
=G\nu_j, \quad j=1,\dots,n.
\end{equation}

Since it is essential that the $\delta$-derivative is
computed on a surface, and $S$ remains constant on this surface then
$\frac{\delta S}{\delta t}=0$. Thus we have
$$
0=\frac{\delta S}{\delta t}=\frac{\partial S}{\partial t}
+\sum_{j=1}^n\frac{\delta S}{\delta x_j}\frac{\delta x_j}{\delta t}
=\frac{\partial S}{\partial t}+\sum_{j=1}^nG|\nabla S|\nu_j^2,
$$
i.e.,
\begin{equation}
\label{g-73}
S_t=-G|\nabla S|.
\end{equation}
From this formula we can see that $-G=\frac{S_t}{|\nabla S|}$ can be
interpreted as the time component of the normal vector.

The space-time unit normal to the surface $\Gamma$ is given by
${\bf n}=\frac{(\nu,-G)}{\sqrt{1+G^2}}$, where
$\sqrt{1+G^2}=\frac{|\nabla_{(x,t)} S|}{|\nabla S|}$, \
$\nabla_{(x,t)}=\big(\nabla,\frac{\partial}{\partial t}\big)$.

If $f(x,t)$ is a function defined only on $\Gamma$, its first
order $\delta$-derivatives with respect to the time and space variables
are defined by the following formulas~\cite[5.2.(15),(16)]{Kan}:
\begin{equation}
\label{g-74}
\frac{\delta f}{\delta t}\stackrel{def}{=}\frac{\partial \widetilde{f}}{\partial t}
+G\frac{\partial\widetilde{f}}{\partial\nu},
\qquad
\frac{\delta f}{\delta x_j}\stackrel{def}{=}\frac{\partial\widetilde{f}}{\partial x_j}
-\nu_j\frac{\partial\widetilde{f}}{\partial\nu},
\quad j=1,\dots,n,
\end{equation}
where $\widetilde{f}$ is a smooth extension of $f$ to a neighborhood
of $\Gamma$ in $\bR^n\times\bR$, \ $j=1,\dots,n$, and
$\frac{\partial\widetilde{f}}{\partial\nu}=\nu\cdot\nabla\widetilde{f}$
is a normal derivative.
Thus the gradient tangent to the surface $\Gamma_t$ is defined as
\begin{equation}
\label{g-74-3}
\nabla_{\Gamma_t}=\nabla-\nabla_{\nu}
=\Big(\frac{\delta}{\delta x_1},\dots,\frac{\delta}{\delta x_n}\Big),
\end{equation}
where $\nabla_{\nu}=\nu\big(\nu\cdot\nabla\big)$ is the gradient
along the normal direction to the surface $\Gamma_t$.

Note that the $\delta$-derivatives (\ref{g-74}) depend only on the
values of $f$ on $\Gamma$, i.e., if $f=0$ on $\Gamma$ then
$\frac{\delta f}{\delta x_j}$ and $\frac{\delta f}{\delta t}$ on
$\Gamma$, \ $j=1,\dots,n$.
Indeed, let $(x_0,t_0)\in \Gamma$. If $\nabla_{(x,t)}f(x_0,t_0)=0$
then $\nabla_{\Gamma_t}f(x_0,t_0)=0$ and $\frac{\delta f}{\delta t}(x_0,t_0)=0$,
where $\nabla_{(x,t)}=\big(\nabla,\frac{\partial}{\partial t}\big)$.
If $\nabla f(x_0,t_0)\ne 0$ then in a neighborhood of the point $(x_0,t_0)$
the surface $\Gamma_{t}$ has the unit space normal
$\nu=\frac{\nabla f}{|\nabla f|}$ and
$G=-\frac{\frac{\partial f}{\partial t}}{|\nabla f|}$.
Consequently, $\nabla_{\Gamma_t}f(x_0,t_0)=0$ and
$\frac{\delta f}{\delta t}(x_0,t_0)=0$.
In the sequel we shall drop tilde from $f$.

For a vector $A(x,t)=(A_1(x,t),\dots,A_n(x,t))$ defined only
on $\Gamma_t$, we introduce the {\it surface {\rm(}tangent{\rm)} divergence}
by the following formula
$$
{\rm div}_{\Gamma_t}A=\nabla_{\Gamma_t}\cdot A
=\sum_{j=1}^n\frac{\delta A_j}{\delta x_j}.
$$
The {\it mean curvature} of the surface $\Gamma_t$ is defined as
\begin{equation}
\label{g-84.3}
{\cK}\stackrel{def}{=}-\frac{1}{2}\nabla_{\Gamma_t}\cdot\nu
=-\frac{1}{2}\sum_{j=1}^n\frac{\delta \nu_j}{\delta x_j}
=-\frac{1}{2}\nabla\cdot\nu.
\end{equation}

\subsection{Distributions defined on a surface.}\label{s5.2}
Consider some facts concerning distributions defined on a
surface~\cite[5.2.]{Kan},~\cite[ch.III,\S 1.]{Gel-S},~\cite{Be1},~\cite{Be2}.
The Heaviside function $H(S)$ is introduced by the following definition:
$$
\big\langle H(S), \ \varphi(x,t)\big\rangle=\int_{S \ge 0}\varphi(x,t)\,dx\,dt,
\quad \forall \, \varphi \in {\cD}(\bR^n\times\bR).
$$
According to~\cite[5.3.(1),(2)]{Kan}, we now introduce the delta function
$\delta(S)$ on the surface $\Gamma$, whose action on a test function
$\varphi(x,t)\in {\cD}(\bR^n\times\bR)$ is given by
\begin{equation}
\label{g-106}
\big\langle \delta(S), \ \varphi(x,t) \big\rangle
=\int_{-\infty}^{\infty}\int_{\Gamma_t}\varphi(x,t)\,d\mu(x)\,dt
=\int_{\Gamma}\varphi(x,t)\frac{\,d\mu(x,t)}{\sqrt{1+G^2}},
\end{equation}
where $d\mu$ is the surface measure on the corresponding surface.
According to~\cite[5.5.Theorem~1.]{Kan}, we have
$$
\frac{\partial H(S)}{\partial x_j}=\nu_j\delta(S),
\qquad
\frac{\partial H(S)}{\partial t}=-G\delta(S).
$$

Now we introduce the derivative of the delta function
$\partial_{\nu}\delta(S)$ along the space normal $\nu$
by the formula~\cite[5.3.(7)]{Kan}
\begin{equation}
\label{g-83}
\bigl\langle \partial_{\nu}\delta(S), \ \varphi \bigr\rangle
=-\Bigl\langle\delta(S), \ \frac{\partial\varphi}{\partial\nu}\Bigr\rangle
=-\int_{-\infty}^{\infty}\int_{\Gamma_t}\frac{\partial\varphi}{\partial\nu}\,d\mu(x)\,dt,
\quad \forall \, \varphi \in {\cD}(\bR^n\times\bR),
\end{equation}
where $\frac{\partial\varphi}{\partial\nu}=\nu\cdot\nabla\varphi$
is the normal derivative of $\varphi$.
If $f(x,t)$ is a continuous function defined on $\Gamma$ which is a
a restriction of some continuous function defined in a neighborhood
of $\Gamma$ in $\bR^n\times\bR$, then the distribution
$\partial_{\nu}(f\delta(S))$ (the so-called {\it double layer}) is a functional
acting by the rule
$$
\bigl\langle \partial_{\nu}\big(f\delta(S)\big), \ \varphi \bigr\rangle
=-\Bigl\langle\delta(S), \ f\frac{\partial\varphi}{\partial\nu}\Bigr\rangle,
\quad \forall \, \varphi \in {\cD}(\bR^n\times\bR).
$$

According to~\cite[5.3.(6)]{Kan}, we have
$$
\delta'(S)=\sum_{i=1}^n\nu_i\frac{\partial}{\partial x_i}\delta(S)
=2{\cK}\delta(S)+\partial_{\nu}\delta(S)
$$
and
\begin{equation}
\label{g-84*}
\frac{\partial}{\partial t}\delta(S)=
-G\big(2{\cK} \delta(S)+\partial_{\nu}\delta(S)\big),
\qquad
\frac{\partial}{\partial x_j}\delta(S)=
\nu_j\big(2{\cK} \delta(S)+\partial_{\nu}\delta(S)\big),
\end{equation}
where ${\cK}$ is the mean curvature (\ref{g-84.3}) of the surface $\Gamma_t$.

If $f(x,t)$ is a differentiable function, using (\ref{g-74}), (\ref{g-84*}),
one can prove the following relations~\cite[12.6.(15),(16)]{Kan}
\begin{equation}
\label{g-85}
\frac{\partial}{\partial x_j}\big(f\delta(S)\big)
=\Big(\frac{\partial f}{\partial x_j}-\nu_j\frac{\partial f}{\partial\nu}
+2{\cK} \nu_jf\Big)\delta(S)+\nu_jf\partial_{\nu}\delta(S),
\quad j=1,\dots,n,
\end{equation}
\begin{equation}
\label{g-86}
\frac{\partial}{\partial t}\big(f\delta(S)\big)
=\Big(\frac{\partial f}{\partial t}+G\frac{\partial f}{\partial\nu}
-2{\cK} Gf\Big)\delta(S)-Gf\partial_{\nu}\delta(S).
\qquad\qquad\qquad\quad
\end{equation}

\subsection{One integrating-by-parts formula.}\label{s5.3}

\begin{Lemma}
\label{g-lem4}
Suppose that $e(x,t)$ is a compactly supported smooth function defined
only on the surface $\Gamma=\bigl\{(x,t):S(x,t)=0\bigr\}$, and $e(x,t)$
is the restriction of some smooth function defined in a neighborhood
of $\Gamma$ in $\bR^n\times\bR$, and $\Gamma_0=\bigl\{x: S(x,0)=0\bigr\}$.
Then the following formula for integration by parts holds:
\begin{equation}
\label{g-84.20}
\int_{\Gamma}e \frac{\delta\varphi}{\delta t}\frac{\,d\mu(x,t)}{\sqrt{1+G^2}}=
\displaystyle
-\int_{\Gamma}\frac{\delta^* e}{\delta t}\varphi\frac{\,d\mu(x,t)}{\sqrt{1+G^2}}
-\int_{\Gamma_0}e(x,0)\varphi(x,0)\,d\mu(x),
\end{equation}
for any $\varphi \in {\cD}(\bR^n\times [0, \infty))$,
where $\frac{\delta^*}{\delta t}$ is the adjoint operator defined as
\begin{equation}
\label{g-84.20-1}
\frac{\delta^* e}{\delta t}=\frac{\delta e}{\delta t}-2{\cK} Ge
=\frac{\delta e}{\delta t}+\nabla_{\Gamma_t}\cdot(eG\nu),
\end{equation}
${\cK}$ is the mean curvature {\rm(\ref{g-84.3})} of the surface $\Gamma_t$.
\end{Lemma}

\begin{proof}
With the help of formulas (\ref{g-106}), (\ref{g-83}), (\ref{g-84*}), (\ref{g-85}),
(\ref{g-86}), we derive by simple calculations
$$
\int_{\Gamma}e\frac{\delta \varphi}{\delta t}\frac{\,d\mu(x,t)}{\sqrt{1+G^2}}
=\Bigl\langle e\delta(S), \ \frac{\delta \varphi}{\delta t}\Bigr\rangle
=\Bigl\langle e\delta(S)H(t), \
\frac{\partial \varphi}{\partial t}+G\frac{\partial\varphi}{\partial\nu}\Bigr\rangle
\qquad\qquad\qquad\qquad
$$
$$
\qquad
=-\Bigl\langle \frac{\partial}{\partial t}\big(e\delta(S)H(t)\big),
\ \varphi \Bigr\rangle
-\Bigl\langle \partial_{\nu}\big(Ge\delta(S)\big)H(t), \ \varphi\Bigr\rangle
$$
$$
\qquad\qquad\qquad
=-\Bigl\langle \frac{\delta e}{\delta t}\delta(S)
-eG\big(2{\cK} \delta(S)+\partial_{\nu}\delta(S)\big), \ \varphi \Bigr\rangle
-\Bigl\langle e\delta(S)\delta(t), \ \varphi \Bigr\rangle
$$
$$
\qquad
-\Bigl\langle \delta(S)\sum_{k=1}^n\frac{\delta (Ge)}{\delta x_k}\nu_k
+eG\partial_{\nu}\delta(S), \ \varphi\Bigr\rangle
$$
$$
\qquad\qquad\qquad
=-\Bigl\langle \Big(\frac{\delta e}{\delta t}-2{\cK} Ge \Big)\delta(S),
\ \varphi \Bigr\rangle
-\Bigl\langle e(x,0)\delta(S(x,0)), \ \varphi(x,0) \Bigr\rangle,
$$
where $H(t)$ is the Heaviside function. Here we use the obvious relation
$$
\sum_{k=1}^n\frac{\delta (Ge)}{\delta x_k}\nu_k=0.
$$
Using the last relation and formula (\ref{g-84.3}), we calculate
$$
\frac{\delta e}{\delta t}-2{\cK} Ge
=\frac{\delta e}{\delta t}+\sum_{j=1}^n\frac{\delta \nu_j}{\delta x_j}Ge
=\frac{\delta e}{\delta t}+\sum_{j=1}^n\frac{\delta (eG\nu_j)}{\delta x_j}.
$$
\end{proof}

\subsection{Transport theorems.}\label{s5.4}
Here we give the following {\it transport theorems\/}.

\begin{Theorem}
\label{g-th8}
{\rm(\cite[12.8.(3)]{Kan},~\cite{Be1},~\cite{Be3},~\cite{Be4})}
Let $f(x,t)$ be a sufficiently smooth function defined in a moving solid
$\Omega_t$, and let a moving hypersurface $\partial \Omega_t$ be its boundary.
Let $\nu$ be the outward unit space normal to the surface $\partial \Omega_t$
and $W(x,t)$ be the velocity of the point $x$ in $\Omega_t$. Then the volume
transport theorem holds:
$$
\frac{d}{dt}\int_{\Omega_t}f(x,t)\,dx
=\int_{\Omega_t}\frac{\partial f}{\partial t}\,dx
+\int_{\partial \Omega_t}f W\cdot\nu \,d\mu(x)
\qquad\quad
$$
\begin{equation}
\label{g-84.1*}
\qquad
=\int_{\Omega_t}\Big(\frac{\partial f}{\partial t}+{\rm div}(fW)\Big)\,dx.
\end{equation}
\end{Theorem}

\begin{Theorem}
\label{g-th4-tr}
{\rm(~\cite[12.8.(9)]{Kan})}
If $e(x,t)$ is a smooth function defined only on the moving surface
$\Gamma_t=\bigl\{x: S(x,t)=0\bigr\}$ {\rm(}which is the restriction of
some smooth function defined in a neighborhood of $\Gamma_t${\rm)},
then the surface transport theorem holds:
$$
\frac{d}{dt}\int_{\Gamma_t}e(x,t)\,d\mu(x)
\qquad\qquad\qquad\qquad\qquad\qquad\qquad\qquad\qquad\qquad
$$
\begin{equation}
\label{g-84.20-tr}
=\int_{\Gamma_t}\Big(\frac{\delta e}{\delta t}-2{\cK} Ge\Big)\,d\mu(x)
=\int_{\Gamma_t}\Big(\frac{\delta e}{\delta t}
+\nabla_{\Gamma_t}\cdot(eU_{\delta})\Big)\,d\mu(x),
\end{equation}
where $U_{\delta}$ is the velocity of $\Gamma_{t}$.
\end{Theorem}

\begin{proof}
Since according to definition (\ref{g-106}),
$$
m(t)=\int_{\Gamma_t}e(x,t)\,d\mu(x)
=\bigl\langle e(x,t)\delta(S), \ 1\bigr\rangle_{x},
$$
using (\ref{g-86}), we obtain
$$
\dot m(t)=\Bigl\langle \frac{\partial}{\partial}\big(e(x,t)\delta(S)\big), \ 1\Bigr\rangle_{x}
=\Bigl\langle
\Big(\frac{\delta e}{\delta t}-2{\cK} Ge\Big)\delta(S)
-Ge\partial_{\nu}\delta(S), \ 1\Bigr\rangle_{x}
$$
$$
=\Bigl\langle
\Big(\frac{\delta e}{\delta t}-2{\cK} Ge\Big)\delta(S), \ 1\Bigr\rangle_{x}
=\int_{\Gamma_t}\Big(\frac{\delta e}{\delta t}-2{\cK} Ge\Big)\,d\mu(x).
\quad
$$
To complete the proof of the theorem, it remains to use formulas
(\ref{g-84.20-1}) and (\ref{g-4.0-10*}).
\end{proof}

\begin{center}
{\bf Acknowledgements }
\end{center}

The author is greatly indebted to E.Yu.~Panov, V.I.~Polischook,
O.S.~Rozanova for fruitful discussions.

\end{document}